\documentclass[11pt,fleqn]{article}
\usepackage{amsmath,amssymb,graphicx}

\begin{document}
\sloppy
\newtheorem{axiom}{Axiom}[section]
\newtheorem{conjecture}[axiom]{Conjecture}
\newtheorem{corollary}[axiom]{Corollary}
\newtheorem{definition}[axiom]{Definition}
\newtheorem{example}[axiom]{Example}
\newtheorem{lemma}[axiom]{Lemma}
\newtheorem{observation}[axiom]{Observation}
\newtheorem{proposition}[axiom]{Proposition}
\newtheorem{theorem}[axiom]{Theorem}

\newcommand{\proof}{\emph{Proof.}\ \ }
\newcommand{\qed}{~~$\Box$}
\newcommand{\rz}{{\mathbb{R}}}
\newcommand{\nz}{{\mathbb{N}}}
\newcommand{\zz}{{\mathbb{Z}}}
\newcommand{\eps}{\varepsilon}
\newcommand{\cei}[1]{\lceil #1\rceil}
\newcommand{\flo}[1]{\left\lfloor #1\right\rfloor}
\newcommand{\seq}[1]{\langle #1\rangle}

\newcommand{\wqap}{Wiener Max-QAP}
\newcommand{\wtree}{{\sc MaxWiener-Tree}}
\newcommand{\aaa}{\alpha}
\newcommand{\bbb}{\beta}
\newcommand{\ddd}{\delta}

\title{{\bf The Wiener maximum quadratic assignment problem}}
\author{
\sc Eranda \c{C}ela\thanks{Corresponding author, {\tt
    cela@opt.math.tu-graz.ac.at}, Phone: {+43 316 873 5366}, Fax: {\tt +43 316
    873 105366},
Institut f\"ur Optimierung und Diskrete Mathematik, TU Graz, Steyrergasse 30,
A-8010 Graz, Austria}
\and\sc Nina S.\ Schmuck\thanks{{\tt nina.schmuck@tugraz.at}.  
Institut f\"ur Optimierung und Diskrete Mathematik, TU Graz, Steyrergasse 30, A-8010 Graz, Austria}
\and\sc Shmuel Wimer\thanks{{\tt wimers@macs.biu.ac.il}. 
School of Engineering, Bar-Ilan University, Ramat-Gan 52900, Israel}
\and\sc Gerhard J.\ Woeginger\thanks{{\tt gwoegi@win.tue.nl}. 
Department of Mathematics and Computer Science, TU Eindhoven, P.O.\ Box 513, 
5600 MB Eindhoven, Netherlands} 
}
\date{}
\maketitle

\begin{abstract}
We investigate a special case of the maximum quadratic assignment problem where
one matrix is a product matrix and the other matrix is the distance matrix of
a one-dimensional point set.
We show that this special case, which we call the Wiener maximum quadratic assignment 
problem, is NP-hard in the ordinary sense and solvable in pseudo-polynomial time.

Our approach also yields a polynomial time solution for the following problem from
chemical graph theory:
Find a tree that maximizes the Wiener index among all trees with a prescribed degree 
sequence.
This settles an open problem from the literature.

\end{abstract}

\medskip
\section{Introduction}
The \emph{Quadratic Assignment Problem} (QAP) in Koopmans-Beckmann form \cite{KoBe1957}
takes as input two $n\times n$ square matrices $A=(a_{ij})$ and $B=(b_{ij})$ with real
entries, and asks to find a permutation $\pi$ that minimizes (or maximizes) the
objective function
\begin{equation}
\label{eq:qap}
Z_{\pi}(A,B) ~:=~ \sum_{i=1}^n\sum_{j=1}^n~ a_{\pi(i)\pi(j)} \, b_{ij}.
\end{equation}
Here $\pi$ ranges over the set $S_n$ of all permutations of $\{1,2,\ldots,n\}$.
The QAP is a hard and well-studied problem in combinatorial optimization; we refer 
the reader to the book \cite{Cela} by \c{C}ela and the  recent book by
Burkard et al.~\cite{Bur09} for more information on this problem.
One branch of research on the QAP concentrates on the algorithmic behavior of 
strongly-structured special cases; see for instance Burkard \& al \cite{BCRW1998}
or Deineko \& Woeginger \cite{DeWo1998} for typical results in this direction.
In this paper, we will contribute a new topic to this research branch.

\subsection{The Wiener Max-QAP}
\label{ssec:intro-qap}
This special case of the QAP restricts matrix $A$ to be a symmetric
\emph{product matrix}, which means that there are non-negative integers
$\aaa_1\le\aaa_2\le\cdots\le\aaa_n$ such that
\begin{equation}
\label{eq:matrixA}
a_{ij} ~=~ \aaa_i\,\aaa_j \mbox{\qquad for $1\le i,j\le n$.}
\end{equation}
The second matrix in the {\wqap} is the distance matrix of a one-dimensional point set,
which means that there are integers $\bbb_1\le\bbb_2\le\cdots\le\bbb_n$ such that
\begin{equation}
\label{eq:matrixB}
b_{ij} ~=~ |\bbb_i-\bbb_j| \mbox{\qquad for $1\le i,j\le n$.}
\end{equation}
Throughout this paper, a matrix $B$ of the form (\ref{eq:matrixB}) will be called
a \emph{1D-distance matrix}.
The goal in the {\wqap} is to \emph{maximize} the objective value in
(\ref{eq:qap}), which can be rewritten as
\begin{equation}
\label{eq:wqap}
Z_{\pi}(A,B) 
~=~ \sum_{i=1}^n\sum_{j=1}^n~ \aaa_{\pi(i)}\, \aaa_{\pi(j)} \, |\bbb_i-\bbb_j|.
\end{equation}

In this paper, we fully determine the computational complexity of the {\wqap}.
Our results are as follows.
As a negative result, we will prove in Section~\ref{sec:complexity} that the {\wqap} 
is NP-hard in the ordinary sense.
On the positive side, in Section~\ref{sec:structure} we will derive a useful decomposition 
property, and we will show that there always exists an optimal permutation that is
V-shaped.
These positive results are then applied in Section~\ref{sec:algorithm} to derive a
pseudo-polynomial time algorithm for the {\wqap}.

\subsection{The Wiener index of a tree}
\label{ssec:intro-tree}
The \emph{Wiener index} $W(G)$ of a connected undirected graph $G=(V,E)$ is the sum 
of the distances between all pairs of vertices in $V$.
The Wiener index was introduced in 1947 by Harold Wiener \cite{Wiener1947} to 
characterize certain molecular structure properties of saturated hydrocarbons.
We refer the reader to the survey article \cite{DoEnGu2001} of Dobrynin, Entringer \& 
Gutman for comprehensive information on this fundamental graph parameter.
Chemists are often interested in the Wiener index of certain trees, where the vertices
represent atoms and where the vertex degrees correspond to the valencies of the atoms.
Entringer, Jackson \& Snyder \cite{EnJaSn1976} show that among all $r$-vertex trees,
the path $P_r$ has the largest and the star $K_{1,r}$ has the smallest Wiener index.
Fischermann \& al \cite{FHRSV2002} characterize the trees that minimize the Wiener index 
among all trees with $r$ vertices and maximum degree at most $\Delta$, and they also
provide several results on the corresponding maximization question.

Wang \cite{Wang2008} describes a simple greedy algorithm for finding a tree that
minimizes the Wiener index among all trees with a prescribed degree sequence. 
Zhang \& al \cite{ZhXiXuPa2008} derive the same result independently and by different
techniques.
The corresponding maximization question remains open (note that the paper \cite{Wang2008} 
claims a solution to it, and that the corrigendum \cite{Wang2008b} points out a crucial 
mistake that invalidates these claims).
Wang \cite{Wang2008b} writes about the maximization question:
\emph{``While the extremal trees seem to be difficult to find and are not unique, 
an algorithm to find at least one of such trees may exist and may be easier to find.''}

In this paper, we resolve the computational complexity of the maximization question:
there exists a polynomial time algorithm that finds a tree that maximizes the Wiener 
index among all trees with a prescribed degree sequence.
More precisely, we will show in Section~\ref{sec:windex} that this maximization question 
can be modeled as a special case of the {\wqap} (modulo certain minor modifications).
Consequently, the machinery developed in Sections~\ref{sec:structure} and
\ref{sec:algorithm} can be applied to it.
The pseudo-polynomial time complexity of our algorithm turns into a polynomial time 
complexity, since all the involved numbers are moderately small.

\medskip
\section{Complexity of the {\wqap}}
\label{sec:complexity}
In this section, we establish NP-hardness of the {\wqap}.
The proof is done by means of a reduction from the following variant of the partition
problem (see Garey \& Johnson \cite{GaJo1979}) which is well-known to be NP-hard in the
ordinary sense.
\begin{quote}
{\sc Problem:} {\sc Partition}
\\[1.0ex]
{\sc Input:}
A sequence $q_1,\ldots,q_{2k}$ of $2k$ positive integers with $\sum_{i=1}^{2k}q_i=2Q$.
\\[1.0ex]
{\sc Question:}
Does there exist $I\subset\{1,\ldots,2k\}$ with $|I|=k$ and
$\sum_{i\in I}q_i=Q$?
\end{quote}
We construct an instance of the {\wqap} of dimension $n=2k$.
The $n\times n$ product matrix $A$ is defined by $\aaa_i= q_i$, $\forall i=1,2,\ldots,n$, and hence
satisfies $a_{ij}=q_iq_j$.
The 1D-distance matrix $B$ uses the points $\bbb_i=1$ for $1\le i\le k$
and $\bbb_i=2$ for $k+1\le i\le2k$.
Note that $b_{ij}=|\bbb_i-\bbb_j|$ equals 1 if one of the indices $i,j$ lies in the
range $1,\ldots,k$ whereas the other index lies in $k+1,\ldots,2k$; in all other cases
$b_{ij}=0$.

Consider a permutation $\pi\in S_{2k}$, and let $J$ denote the set of all $i$ with
$\pi(i)\le k$.
Note that $|J|=k$.
By setting $x=\sum_{i\in J}q_i$, the objective value in (\ref{eq:wqap}) can be
rewritten as
\[  \sum_{i\in J}\sum_{j\notin J} q_iq_j
~=~ \sum_{i\in J}q_i \sum_{j\notin J}q_j
~=~ x(2Q-x)
~=~ 2Qx - x^2.\]
As the concave function $f(x)=2Qx-x^2$ is maximized at $x=Q$, it is easily seen
that the QAP has objective value at least $Q^2$ if and only if the {\sc Partition}
instance has a positive answer.
This yields the following result.

\begin{theorem}
\label{th:NP-hard}
The {\wqap} is NP-hard in the ordinary sense.
\qed
\end{theorem}

\medskip
\section{Structure of the {\wqap}}
\label{sec:structure}
We first discuss a useful decomposition property of the {\wqap}.
Consider some fixed instance with product matrix $A$ and 1D-distance matrix $B$,
and let $I\subset\{1,\ldots,n\}$ with $|I|=k$ be some fixed subset of the indices.
Instead of optimizing over all possible permutations in (\ref{eq:wqap}), we only
allow permutations $\pi$ that map $I$ into $J=\{1,\ldots,k\}$ and that consequently
map the complement of $I$ into the complement $\{k+1,\ldots,n\}$ of $J$.
Intuitively speaking, this subdivided version QAP assigns the values $\aaa_i$ with
$i\in I$ to the points $\bbb_1,\ldots,\bbb_k$ and the values $\aaa_i$ with $i\notin I$
to the remaining points $\bbb_{k+1},\ldots,\bbb_n$.

Any permutation $\pi\in S_n$ induces a bijection $\sigma$ that maps $J$ into $I$,
and another bijection $\tau$ that maps the complement of $J$ into the complement
of $I$.
Furthermore, denote $X=\sum_{i\in I}\aaa_i$ and $Y=\sum_{i\notin I}\aaa_i$.
The objective value $Z_{\pi}(A,B)$ in (\ref{eq:wqap}) can then be written as 
$Z_1+Z_2+Z_3+Z_4$ in the following way.

\begin{eqnarray}
\label{eq:a.Z1}
Z_1
&=& \sum_{i\in J}\sum_{j\in J} ~ \aaa_{\pi(i)}\aaa_{\pi(j)} \, |\bbb_i-\bbb_j|
~=~ \sum_{i\in J}\sum_{j\in J} ~ \aaa_{\sigma(i)}\aaa_{\sigma(j)} \, |\bbb_i-\bbb_j|
\end{eqnarray}

\begin{eqnarray}
\label{eq:a.Z2}
Z_2
&=& \sum_{i\notin J}\sum_{j\notin J} ~ \aaa_{\pi(i)}\aaa_{\pi(j)} \, |\bbb_i-\bbb_j|
~=~ \sum_{i\notin J}\sum_{j\notin J} ~ \aaa_{\tau(i)}\aaa_{\tau(j)} \, |\bbb_i-\bbb_j|
\end{eqnarray}

\begin{eqnarray}
Z_3
&=& \sum_{i\in J}\sum_{j\notin J} ~ \aaa_{\pi(i)}\aaa_{\pi(j)} \, |\bbb_i-\bbb_j|
~=~ \sum_{i\in J}\sum_{j\notin J} ~ \aaa_{\sigma(i)}\aaa_{\tau(j)} \, (\bbb_j-\bbb_k+\bbb_k-\bbb_i)
\nonumber
\\[0.5ex]
&=&
\sum_{i\in J} ~ \aaa_{\sigma(i)} Y (\bbb_k-\bbb_i) ~+~
\sum_{j\notin J} ~\aaa_{\tau(j)} X (\bbb_j-\bbb_k)
\label{eq:a.Z3}
\end{eqnarray}

\begin{eqnarray}
Z_4
&=& \sum_{i\notin J}\sum_{j\in J} ~ \aaa_{\pi(i)}\aaa_{\pi(j)} \, |\bbb_i-\bbb_j|
~=~ Z_3
\label{eq:a.Z4}
\end{eqnarray}

Note that in the resulting summations in (\ref{eq:a.Z1})--(\ref{eq:a.Z3}), every single 
term does either depend on function $\sigma$ or on function $\tau$, but does never 
depend on both functions simultaneously.
If we collect all the terms in $Z_1+Z_2+Z_3+Z_4$ that solely depend on this function
$\sigma$, we get
\begin{equation}
\label{eq:a.1}
\sum_{i\in J}\sum_{j\in J} ~ \aaa_{\sigma(i)}\aaa_{\sigma(j)} \, |\bbb_i-\bbb_j|
+ 2 \sum_{i\in J} ~ \aaa_{\sigma(i)} Y (\bbb_k-\bbb_i).
\end{equation}
We observe that the objective function in (\ref{eq:a.1}) essentially corresponds to a
smaller $(k+1)$-dimensional {\wqap}:
the underlying 1D-distance matrix is built around the $k$ points $\bbb_1,\ldots,\bbb_k$
plus a duplicated point at $\bbb_k$ (which corresponds to the occurrence of $\bbb_k$ in
the right hand sum).
The underlying product matrix is built around the $k$ numbers $\aaa_i$ with $i\in I$
plus the number $Y$ (which corresponds to the factor $Y$ in the right hand sum).
Furthermore, the right hand sum in (\ref{eq:a.1}) imposes the additional restriction 
that the new value $Y$ has to be assigned to the duplicated point at $\bbb_k$.

As a consequence of all this, the problem of finding the optimal function $\sigma$ and
the problem of finding the optimal function $\tau$ are two separate optimization problems
that can be solved independently of each other.
We call this the \emph{decomposition property} of the {\wqap}.
This decomposition property plays a central role in many of our arguments.
As a first application of the decomposition property, we next deduce a result on the
combinatorial structure of optimal permutations for the {\wqap}.

\begin{definition}
\label{df:v-shaped}
A permutation $\pi\in S_n$ is called \emph{V-shaped}, if there exists an index $\ell$
with $1\le\ell\le n$ such that $\pi(i)>\pi(i+1)$ for $i=1,\ldots,\ell-1$ and
such that $\pi(i)<\pi(i+1)$ for $i=\ell,\ldots,n-1$.
\end{definition}
In other words, a V-shaped permutation $\pi$ is first decreasing up to $\ell$, and then
increasing from $\ell$ onwards, where the increasing or the decreasing part can
also be empty. 

\begin{theorem}
\label{th:v-shaped}
Every instance of the {\wqap} possesses an optimal solution $\pi$ that is V-shaped.
\end{theorem}
\proof
For simplicity of presentation, we will assume without much loss of generality that all
$\aaa$-values are pairwise distinct, and that also all $\bbb$-values are pairwise distinct.
The statement for the general case then follows easily from this (by locally reordering 
or renaming the values).

Now consider an optimal permutation $\pi$, and note that the values $\pi(1),\ldots,\pi(n)$
and $\aaa_{\pi(1)},\ldots,\aaa_{\pi(n)}$ are ordered in the same way.
Suppose for the sake of contradiction that permutation $\pi$ has a local maximum at $k$ 
with $\aaa_{\pi(k-1)}<\aaa_{\pi(k)}$ and $\aaa_{\pi(k+1)}<\aaa_{\pi(k)}$.
By the decomposition property of the {\wqap}, the optimal permutation $\pi$ induces
an optimal solution to the five-dimensional problem of assigning the five values
\[    L=\sum_{i=1}^{k-2}\aaa_{\pi(i)},
\quad \aaa_{\pi(k-1)},
\quad \aaa_{\pi(k)},
\quad \aaa_{\pi(k+1)},
\quad R=\sum_{i=k+2}^{n}\aaa_{\pi(i)}
\]
to the five points $\bbb_{k-1}$, $\bbb_{k-1}$, $\bbb_k$, $\bbb_{k+1}$, and $\bbb_{k+1}$
subject to the constraint that value $L$ is assigned to point $\bbb_{k-1}$ and that value
$R$ is assigned to point $\bbb_{k+1}$.

Now let us switch the positions of $\aaa_{\pi(k-1)}$ and $\aaa_{\pi(k)}$ in the solution
that $\pi$ induces for the five-dimensional problem, such that $\aaa_{\pi(k-1)}$ goes to
$\bbb_k$ and $\aaa_{\pi(k)}$ goes to $\bbb_{k-1}$.
Since this switch cannot increase the objective value, the difference between the corresponding
two objective values is non-positive:
\begin{equation}
\label{eq:b.1}
\left(\aaa_{\pi(k)}-\aaa_{\pi(k-1)}\right) \,
\left(\bbb_k-\bbb_{k-1}\right) \,
\left(\aaa_{\pi(k+1)}+R-L\right) ~\le~ 0
\end{equation}
In an analogous way, we can switch the positions of $\aaa_{\pi(k)}$ and $\aaa_{\pi(k+1)}$
in the induced solution.
This then leads to the following inequality:
\begin{equation}
\label{eq:b.2}
\left(\aaa_{\pi(k)}-\aaa_{\pi(k+1)}\right) \,
\left(\bbb_{k+1}-\bbb_k\right) \,
\left(\aaa_{\pi(k-1)}-R+L\right) ~\le~ 0
\end{equation}
The first two factors on the left-hand side of (\ref{eq:b.1}) and also the first two
factors on the left-hand side of (\ref{eq:b.2}) all are positive.
This implies $\aaa_{\pi(k+1)}+R-L\le0$ and $\aaa_{\pi(k-1)}-R+L\le0$.
By summing these two inequalities we arrive at the contradiction
$\aaa_{\pi(k+1)}+\aaa_{\pi(k-1)}\le0$.
We conclude that $\pi$ cannot have any local maximum, and this yields that $\pi$ indeed
is V-shaped.
\qed

\medskip
\section{An algorithm for the {\wqap}}
\label{sec:algorithm}
In this section, we design a pseudo-polynomial time algorithm for the {\wqap} that is 
based on a standard dynamic programming approach.
We recall that the two underlying sequences $\aaa_1\le\aaa_2\le\cdots\le\aaa_n$ and 
$\bbb_1\le\bbb_2\le\cdots\le\bbb_n$ are in non-decreasing order.

Every state in the dynamic program is specified by a quadruple $(k,m,L,R)$ of 
integers that satisfy the following conditions:
\begin{equation}
\label{eq:dp.1}
1\le k\le n, \qquad 
1\le m\le n-k+1, \qquad 
0\le L,R, \qquad
L+R=\sum_{i=k+1}^n\aaa_i. 
\end{equation}

With every such state $(k,m,L,R)$ we associate the following $(k+2)$-dimensional {\wqap}:
the product matrix results from the $k+2$ non-negative integers in the sequence 
$\aaa_1,\ldots,\aaa_k,L,R$, 
and the 1D-distance matrix results from the $k$ points $\bbb_m,\bbb_{m+1},\ldots,\bbb_{m+k-1}$ 
plus another point in $\bbb_m$ plus another point in $\bbb_{m+k-1}$.
The goal is to find the best solution to this QAP subject to the constraint that the value
$L$ is assigned to point $\bbb_m$ and that the value $R$ is assigned to $\bbb_{m+k-1}$.
In other words, we want to find a bijection $\sigma$ from $\{m,\ldots,m+k-1\}$ to
$\{1,\ldots,k\}$ that maximizes the objective value 
\begin{eqnarray*}
Z &=&
\sum_{i=m}^{m+k-1} \sum_{j=m}^{m+k-1} \aaa_{\sigma(i)}\aaa_{\sigma(j)} \, |\bbb_i-\bbb_j|
~+~ 2LR\,\, |\bbb_{m+k-1} -\bbb_m| ~+~
\\[1.0ex]
&& ~~~+~
2\sum_{i=m}^{m+k-1} \aaa_{\sigma(i)} L\, |\bbb_i-\bbb_m| ~+~
2\sum_{i=m}^{m+k-1} \aaa_{\sigma(i)} R\, |\bbb_{m+k-1} -\bbb_i| 
\end{eqnarray*}
We use $Z(k,m,L,R)$ to denote the maximum objective value of the corresponding state.
Next, we will describe how to compute and to store all the values $Z(k,m,L,R)$ step 
by step and in increasing order of $k$.
For $k=1$ the corresponding instances are trivial to solve, since they only have a 
single feasible solution.

Next consider some fixed state $(k,m,L,R)$ with $k\ge2$, and denote 
$M=\sum_{i=1}^{k-1}\aaa_i$.
By the decomposition property discussed in Section~\ref{sec:structure} and by 
Theorem~\ref{th:v-shaped}, there exists an optimal bijection $\sigma$ that induces a 
V-shaped assignment of the $k$ values $\aaa_1,\ldots,\aaa_k$ to the $k$ points 
$\bbb_m,\bbb_{m+1},\ldots,\bbb_{m+k-1}$.
This implies that $\aaa_k$ as the largest of the $k$ values must be assigned either 
to point $\bbb_m$ or to point $\bbb_{m+k-1}$.
First consider the case where $\aaa_k$ is assigned to point $\bbb_m$. 
Then the remaining $k-1$ values are assigned to $\bbb_{m+1},\ldots,\bbb_{m+k-1}$, and
by the decomposition property the largest possible objective value in this case is
\begin{equation}
\label{eq:dp.2a}
Z_1 ~:=~ Z(k-1,m+1,L+\aaa_k,R) ~+~ 2(L+\aaa_k)\,(M+R)\,|\bbb_{m+1}-\bbb_m|.
\end{equation}
In the second case the value $\aaa_k$ is assigned to point $\bbb_{m+k-1}$.
Then the remaining $k-1$ values are assigned to $\bbb_m,\ldots,\bbb_{m+k-2}$, and
the largest possible objective value is
\begin{equation}
\label{eq:dp.2b}
Z_2 ~:=~ Z(k-1,m,L,R+\aaa_k) ~+~ 2(L+M)\,(R+\aaa_k)\,|\bbb_{m+k-1}-\bbb_{m+k-2}|.
\end{equation}
This yields $Z(k,m,L,R)=\max\{Z_1,Z_2\}$, and in this fashion one easily determines 
all the function values $Z(k,m,L,R)$ with $2\le k\le n$.
In the end, the optimal objective value of the underlying QAP instance can be found
as $Z(n,1,0,0)$.

\begin{theorem}
\label{th:dp}
The {\wqap} has a pseudo-polynomial time solution algorithm with time complexity
$O(n^2\cdot\sum\aaa_i)$.
\end{theorem}
\proof
The correctness of the dynamic programming approach is clear from the above considerations.
It remains to analyze the time complexity.
We observe that there are only $O(n^2\cdot\sum\aaa_i)$ different states $(k,m,L,R)$
in the dynamic program:
there are $O(n)$ possible values for $k$ and $m$, respectively, and there are 
$O(\sum\aaa_i)$ possible values for $L$; note that the value of $R$ is already fully 
determined by the values of $k$ and $L$.

The expressions in (\ref{eq:dp.2a}) and (\ref{eq:dp.2b}) can be evaluated in constant
time $O(1)$: they refer to values $Z(k-1,*,*,*)$ that are known from earlier stages of
the dynamic program, and they refer to the values $M=\sum_{i=1}^{k-1}\aaa_i$ that can
all be precomputed and stored in a preprocessing phase.
All in all, this yields that the time complexity is proportional to the number of
states and hence is $O(n^2\cdot\sum\aaa_i)$.
\qed

Note that our approach only yields the optimal objective value.
By storing appropriate auxiliary information in the states of the dynamic program, one 
can also compute the corresponding optimal permutation within the same time complexity.
These are standard techniques, and we do not elaborate on them.

\medskip
\section{Maximizing the Wiener index of a tree}
\label{sec:windex}
We now return to the Wiener index of a graph that has been introduced and discussed 
in Section~\ref{ssec:intro-tree}.
We will investigate the following algorithmic problem {\wtree} that was left open 
by Wang \cite{Wang2008,Wang2008b}:  
an instance consists of a degree sequence $d_1,\ldots,d_r$ of $r$ positive integers 
with $\sum_{i=1}^rd_i=2r-2$.
The goal is to determine the largest possible Wiener index over all trees with 
degree sequence $d_1,\ldots,d_r$.

Recall that a \emph{caterpillar} is a tree that turns into a path (the so-called
\emph{backbone} of the caterpillar) if all its leaves are removed.
Shi \cite{Shi1993} proved that for every instance of {\wtree}, all the maximizing 
trees are caterpillars; other proofs of this result can be found for instance in 
Wang \cite{Wang2008} and Schmuck \cite{Schmuck2010}.

Now consider such a caterpillar $T$, let $v_1,\ldots,v_n$ denote the vertices ordered 
along the backbone of the caterpillar, and let $\ell_i$ with $1\le i\le n$ denote the 
number of leaves adjacent to vertex $v_i$.
We define set $C_i$ to consist of vertex $v_i$ together with its $\ell_i$ adjacent leaves.
Then the vertex pairs inside $C_i$ contribute an amount of $\ell_i^2$ to the Wiener 
index $W(T)$, and the pairs with one vertex in $C_i$ and one vertex in $C_j$ 
($i\ne j$) contribute
\begin{eqnarray*}
\lefteqn{\ell_i\ell_j(|j-i|+2) ~+~ \ell_i\,(|j-i|+1) ~+~ \ell_j\,(|j-i|+1) ~+~ |j-i|}
\\[1.0ex] && \qquad=~~
(\ell_i+1)\,(\ell_j+1)\,|j-i| ~+~ (2\ell_i\ell_j+\ell_i+\ell_j).
\end{eqnarray*}
Hence the Wiener index of this caterpillar $T$ is
\begin{eqnarray}
W(T) &=& \sum_{i=1}^n\ell_i^2 ~+~
\sum_{i=1}^n \sum_{j=i+1}^n [ (\ell_i+1)\,(\ell_j+1)\,|j-i| ~+~ (2\ell_i\ell_j+\ell_i+\ell_j)]
\nonumber
\\[1.0ex]
&=& \left(\sum_{i=1}^n\ell_i\right)^2+(n-1)\sum_{i=1}^n\ell_i ~+~
\frac12 \sum_{i=1}^n \sum_{j=1}^n (\ell_i+1)\,(\ell_j+1)\,|j-i|.
\label{eq:t.1}
\end{eqnarray}
Neither the first nor the second sum in (\ref{eq:t.1}) do depend on the way how the 
numbers $\ell_1,\ldots,\ell_n$ are assigned to the backbone vertices $v_1,\ldots,v_n$.
The assignment in the third sum in (\ref{eq:t.1}) yields an instance of the {\wqap}
with values $\aaa_i=\ell_i+1$ and points $\bbb_i=i$.

\begin{observation}
\label{ob:caterpillar}
The problem of finding the maximizing caterpillar for an instance of {\wtree} with an 
explicitly specified backbone $v_1,\ldots,v_n$ and an explicitly specified sequence
$\ell_1,\ldots,\ell_n$ of leaf-numbers is equivalent to a {\wqap}.
\qed
\end{observation}

Next, consider a degree sequence $d_1,\ldots,d_r$ that forms an instance of {\wtree}, 
and assume without loss of generality that $2\le d_1\le d_2\le\cdots\le d_n$ and that 
$d_{n+1}=\cdots=d_r=1$.
(Since the case $n=1$ is trivial, we assume from now on that $n\ge2$.)
It is straightforward to see that the backbone of the maximizing caterpillar will 
consist of $n$ vertices $v_1,\ldots,v_n$.
In contrast to this, it is not straightforward to write down the sequence 
$\ell_1,\ldots,\ell_n$ of leaf-numbers:
If one of the inner backbone vertices $v_i$ with $2\le i\le n-1$ gets degree $d_i$ 
then it is adjacent to $\ell_i=d_i-2$ leaves, whereas
if one of the two outermost backbone vertices $v_i$ with $i=1$ or $i=n$ gets degree 
$d_i$ then it is adjacent to $\ell_i=d_i-1$ leaves.

Motivated by the discussion in the preceding paragraph, we introduce an 
$(n+2)$-dimensional instance of the {\wqap}:
the product matrix is built around the $n$ numbers $\aaa_i=d_i-1$ for $1\le i\le n$
and the two additional numbers $\aaa_{n+1}=\aaa_{n+2}=1$.
The 1D-distance matrix is built around the $n$ points $\bbb_i=i-1$ for $2\le i\le n+1$
and the two additional points $\bbb_1=1$ and $\bbb_{n+2}=n$.
Furthermore, we impose the constraints that value $\aaa_{n+1}=1$ must be assigned to
point $\bbb_1=1$, and that value $\aaa_{n+2}=1$ must be assigned to point $\bbb_{n+2}=n$.
These additionally imposed constraints take care of the special treatment of the two 
outermost backbone vertices $v_1$ and $v_n$.

Now let us verify that the machinery of Sections~\ref{sec:structure} and 
\ref{sec:algorithm} still can be applied to this variant of the {\wqap}.
First of all, the decomposition property works out exactly as before. 
Also Theorem~\ref{th:v-shaped} continues to hold, since adding $\pi(1)=n+1$ and 
$\pi(n+2)=n+2$ to a V-shaped permutation $\pi(2),\ldots,\pi(n+1)$ of the numbers
$1,\ldots,n$ always yields a V-shaped permutation.
The dynamic program in Section~\ref{sec:algorithm} needs some small cosmetic changes
that are caused by the additionally imposed constraints.
\begin{itemize}
\itemsep=-0.3ex
\item
All states $(k,m,L,R)$ with $k\le n$ must satisfy $2\le m\le n-k+2$ and $L,R\ge1$.
\item
All states $(n+1,m,L,R)$ must satisfy $m=1$, $L=0$, and $R=1$.
\item
All states $(n+2,m,L,R)$ must satisfy $m=1$ and $L=R=0$.
\end{itemize}
These conditions ensure that the dynamic program assigns the values $\aaa_{n+1}$ and 
$\aaa_{n+2}$ during the last two stages to points $\bbb_1$ and $\bbb_{n+2}$ exactly 
as desired.

Furthermore, we note that in equations (\ref{eq:dp.2a}) and (\ref{eq:dp.2b}) we always 
have $|\bbb_{m+1}-\bbb_m|=1$ and $|\bbb_{m+k-1}-\bbb_{m+k-2}|=1$ for $k\le n$, and 
$|\bbb_{m+1}-\bbb_m|=0$ and $|\bbb_{m+k-1}-\bbb_{m+k-2}|=|\bbb_k-\bbb_{k-1}|$
 for $k>n$.
This implies that neither the recursive computations nor the values $Z(k,m,L,R)$ 
do depend on the second coordinate $m$, which consequently may be dropped.
(This should also be clear intuitively, since this second coordinate encodes the piece 
of the backbone to which the first $k$ values $\aaa_1,\ldots,\aaa_k$ are assigned. 
All backbone pieces of length $k$ are paths on $k$ vertices and thus have the same 
combinatorial structure.)

What about the time complexity?
Exactly as in the proof of Theorem~\ref{th:dp} the time complexity is proportional 
to the number of states $(k,L,R)$.
Since $k$ can take $O(r)$ possible values, and since $L$ can take $O(\sum\aaa_i)=O(r)$ 
possible values, the time complexity is $O(r^2)$.

\begin{theorem}
\label{th:tree}
The problem of finding a tree that maximizes the Wiener index among all trees with a 
prescribed degree sequence can be solved in quadratic time $O(r^2)$, where $r$ denotes the overall number of terms in the degree sequence.
\qed
\end{theorem}

\medskip
\section{Conclusions}
\label{sec:conclusion}
We have introduced the {\wqap}, a special case of the quadratic assignment problem. 
We have provided a complete picture of the computational complexity of this special case:
It is NP-hard in the ordinary sense, and it is solvable in pseudo-polynomial time.
Our investigations also gave us a polynomial time algorithm for finding a tree that 
maximizes the Wiener index among all trees with a prescribed degree sequence, 
thereby settling a prominent open problem from chemical graph theory.

One obvious open problem is to bring the quadratic time complexity $O(r^2)$ in 
Theorem~\ref{th:tree} down to $O(r\log r)$ or perhaps even down to linear time $O(r)$.

Another open problem concerns the Wiener Min-QAP, where the goal is to \emph{minimize}
the objective value in (\ref{eq:wqap}).
It is an easy exercise to rewrite and to adapt the results of 
Sections~\ref{sec:structure} and \ref{sec:algorithm} to the minimization version:
The minimization version always has an optimal solution that is \emph{pyramidal}
(which means that the permutation is first increasing up to some value $\ell$, 
and then decreasing from $\ell$ onwards).
And the minimization problem can be solved by dynamic programming in pseudo-polynomial 
time, within the same time complexity as that in Theorem~\ref{th:dp}.
The main gap in our knowledge concerns the complexity of the Wiener Min-QAP,
and we pose the open problem of deciding whether it actually is NP-hard.

\section*{Acknowledgements}
\nopagebreak
Gerhard Woeginger acknowledges support
by the Netherlands Organization for Scientific Research (NWO), grant 639.033.403,
by DIAMANT (an NWO mathematics cluster),
and by BSIK grant 03018 (BRICKS: Basic Research in Informatics for Creating the
Know\-ledge Society).
\bigskip

\noindent
We thank Stephan Wagner  (Stellenbosch University) for bringing the
problem of maximizing the Wiener
Index  of a tree to our attention.


\end{document}